\documentclass[12pt]{amsart}
\usepackage{somedefs,tracefnt,latexsym,,rawfonts,epsfig,amsfonts,amssymb,latexsy
m,enumerate,amscd}
\def\cal{\mathcal}
\def\Bbb{\mathbb}

\newtheorem{prop}{Proposition}[section]
\newtheorem{thm}{Theorem}[section]

\newtheorem{lemma}{Lemma}[section]
\newtheorem{cor}{Corollary}[section]

\newtheorem{defn}{Definition}[section]

\newtheorem{rem}{Remark}[section]
\begin{document}
\title[Algebraic $K$-theory] {Algebraic $K$-theory of groups wreath 
product with finite groups}
\author[S.K. Roushon]{S.K. Roushon*}
\date{June 22, 2006. Final revision on January 31, 2007}
\address{School of Mathematics\\
Tata Institute\\
Homi Bhabha Road\\
Mumbai 400005, India}
\email{roushon@math.tifr.res.in} 
\urladdr{http://www.math.tifr.res.in/\~\ roushon/paper.html}
\thanks{*This work was done during the author's visit to the 
Mathematisches Institut, Universit\"{a}t M\"{u}nster under a 
Fellowship of the Alexander von Humboldt-Stiftung}
\begin{abstract} The 
Farrell-Jones Fibered Isomorphism Conjecture for the stable 
topological pseudoisotopy 
theory has been proved for several classes of groups. For example for  
discrete 
subgroups of Lie groups (\cite{FJ}), virtually poly-infinite cyclic 
groups (\cite{FJ}), Artin 
braid groups (\cite{FR}), a class of virtually poly-surface groups 
(\cite{R3}) and virtually solvable 
linear group (\cite{FL}). We extend these results in the 
sense that if $G$ is a group from the above classes then we prove the 
conjecture 
for the wreath product $G\wr H$ for $H$ a finite group. The need for 
this kind of extension is already evident in \cite{FR}, \cite{R2} and 
\cite{R3}. We also prove the conjecture for some other classes of 
groups.\end{abstract}



\maketitle

\section{Introduction} In this article we are mainly concerned about the 
Fibered Isomorphism Conjecture for the stable topological pseudoisotopy
theory. We extend some existing results and also prove the conjecture for
some other classes of groups. Finally we deduce a corollary for the
Isomorphism Conjecture for the algebraic $K$-theory in dimension $\leq 1$
(see Corollary \ref{kth}).  

Before we state the theorem let us recall that given two groups $G$ and 
$H$ the {\it wreath product} 
$G\wr H$ is by definition the semidirect product $G^H\rtimes H$ where the 
action of $H$ on $G^H$ is the regular action. 

\begin{thm} \label{intro}
The Fibered Isomorphism Conjecture for the stable topological 
pseudoisotopy theory is true 
for the group $G\wr H$ where $H$ is a finite group and $G$ is one of the 
following groups. 

(a). Virtually polycyclic groups. 

(b). Virtually solvable subgroups of $GL_n({\Bbb C})$.

(c). Cocompact discrete subgroups of virtually connected Lie groups.

(d). Artin full braid groups.

(e). Weak strongly poly-surface groups (see Definition \ref{defn}).

(f). Extensions of closed surface groups by surface groups.

(g). $\pi_1(M)\rtimes \Bbb Z$, where $M$ is a closed Seifert fibered 
space.\end{thm}

\begin{proof} See Corollary \ref{cons} and Remark 
\ref{remark}.\end{proof}

(In the notation defined in Definition \ref{definition} 
the Theorem says that the 
FIC$^{_P{\cal H}^?_*}_{\cal {VC}}$ is true for $G\wr H$ or equivalently 
the
FICwF$^{_P{\cal H}^?_*}_{\cal {VC}}$ is true for $G$.)

The FICwF$^{_P{\cal H}^?_*}_{\cal {VC}}$ was proved for $3$-manifold 
groups in \cite{R2} and 
\cite{R3} and 
for the fundamental groups of a certain class of graphs of virtually  
poly-cyclic groups in \cite{R4}.

For the extended Fibered Isomorphism Conjecture FICwF$^{_P{\cal 
H}^?_*}_{\cal {VC}}$ we deduce the 
following proposition using $(a)$ of Lemma \ref{inverse}. This result is 
not yet known 
if 
we replace FICwF$^{_P{\cal 
H}^?_*}_{\cal {VC}}$ by FIC$^{_P{\cal H}^?_*}_{\cal {VC}}$ (see 
5.4.5 in 
\cite{LR} for the question).

\begin{prop} \label{corprop} Let $G$ be a group containing a finite index 
subgroup 
$\Gamma$. Assume that the FICwF$^{_P{\cal H}^?_*}_{\cal {VC}}$ is true for 
$\Gamma$. Then the 
FICwF$^{_P{\cal H}^?_*}_{\cal {VC}}$ 
is true for $G$.\end{prop}


We work in the general setting of the conjecture in equivariant 
homology theory (see \cite{BL}). We 
find out a set of properties which are all satisfied in the pseudoisotopy 
case of the conjecture. And assuming these properties we prove a theorem 
(Theorem \ref{mth}) for the Isomorphism Conjecture in equivariant homology 
theory.
Theorem 
\ref{intro} is then a particular case of Theorem \ref{mth}. We mention in 
Corollary \ref{kth} another consequence of our result for the 
Isomorphism Conjecture in the algebraic $K$-theory case. We also hope that 
the general Theorem \ref{mth} will be useful for future application.

The methods we use in this paper were developed in \cite{R4}.

\medskip
\noindent
{\bf Acknowledgement.} The author would like to thank the referee for
carefully reading the paper and for some suggestions. 

\section{Statement of the general Theorem} We first recall 
the general statement of 
the Isomorphism Conjecture in equivariant homology theory from 
\cite{BL} and also recall some definitions from \cite{R4}.

A class $\cal C$ of subgroups of a group $G$ is called a {\it family of 
subgroups} if $\cal C$ is closed under taking subgroups and conjugations.
For a family of subgroups $\cal C$ of $G$, $E_{\cal C}(G)$ denotes a 
$G$-CW complex so that the fixpoint set $E_{\cal C}(G)^H$ is 
contractible if $H\in {\cal C}$ and empty otherwise. And $R$ denotes an 
associative ring with unit. 

\medskip
\noindent
{\bf (Fibered) Isomorphism Conjecture:} (see definition 1.1 in 
\cite{BL}) Let ${\cal H}^?_*$ be an 
equivariant homology theory with values in $R$-modules. Let $G$ be a group 
and $\cal C$ be a family of subgroups of $G$. Then the {\it 
Isomorphism Conjecture} for the pair $(G, {\cal C})$ states that the 
projection 
$p:E_{\cal C}(G)\to pt$ to the point $pt$ induces an isomorphism $${\cal 
H}^G_n(p):{\cal H}^G_n(E_{\cal C}(G))\simeq {\cal H}^G_n(pt)$$ for $n\in 
{\Bbb Z}$. 
And the {\it
Isomorphism Conjecture in dimension $\leq k$}, for $k\in {\Bbb Z}$, states 
that ${\cal H}^G_n(p)$ is an isomorphism for $n\leq k$.
Finally the {\it Fibered Isomorphism Conjecture} for the pair $(G, {\cal 
C})$ states that given a group homomorphism $\phi: K\to G$ the 
Isomorphism Conjecture is true for the pair $(K, \phi^*{\cal C})$.
Here $\phi^*{\cal C}=\{H<K\ |\ \phi(H)\in {\cal C}\}$.
\medskip

From now on let $\cal C$ be a class of groups which is closed under 
isomorphisms, 
taking subgroups and taking quotients. The class $\cal {VC}$ of all 
virtually cyclic groups has these properties. Given a group $G$ we 
denote by 
${\cal C}(G)$ the class of subgroups of $G$ belonging to $\cal C$. Then  
${\cal C}(G)$ is a family of subgroups of $G$ which is also closed under 
taking quotients.

\begin{defn} \label{definition} (see definition 2.1 in \cite{R4}) {\rm If 
the (Fibered) Isomorphism Conjecture is
true for the pair $(G, {\cal C}(G))$ we say that the {\it (F)IC$^{{\cal
H}^?_*}_{\cal C}$ is true for $G$} or simply say {\it (F)IC$^{{\cal
H}^?_*}_{\cal C}(G)$ is satisfied}. Also we say that the {\it 
(F)ICwF$^{{\cal
H}^?_*}_{\cal C}(G)$ is satisfied} if
the {\it (F)IC$^{{\cal
H}^?_*}_{\cal C}$} is true for $G\wr H$ for all finite groups 
$H$.}\end{defn}

Let us denote by $_P{\cal H}^?_*$ and $_K{\cal H}^?_*$ the equivariant 
homology theories arise in the 
Isomorphism Conjecture of Farrell and Jones 
(\cite{FJ}) corresponding to the 
stable topological pseudoisotopy theory and 
the algebraic $K$-theory respectively. 
For these homology theories the 
conjecture is identical with the conjecture made in \S 1.6 and 
\S 1.7 in \cite{FJ}. (See sections 5 and 6 in \cite{BL} for the second  
case and 4.2.1 and 4.2.2 in \cite{LR} for the first case.)

Note that if the FIC$^{{\cal H}^?_*}_{\cal C}$ (respectively 
FICwF$^{{\cal H}^?_*}_{\cal C}$) is true for a group $G$ 
then the FIC$^{{\cal H}^?_*}_{\cal C}$ (respectively
FICwF$^{{\cal H}^?_*}_{\cal C}$) is true for  
subgroups of $G$. We refer to this fact as the {\it hereditary property}.
Also note that the (F)IC$^{{\cal H}^?_*}_{\cal C}$ is true for $H\in 
{\cal C}$.
 
\begin{defn} \label{property} (see definition 2.2 in \cite{R4}) {\rm 
We say that the {\it ${\cal P}^{{\cal H}^?_*}_{\cal C}$-property} is 
satisfied if for $G_1, G_2\in {\cal C}$ the product
$G_1\times G_2$ satisfies the FIC$^{{\cal H}^?_*}_{\cal C}$.}\end{defn}

In the following the notation $A\rtimes B$ stands for the semidirect 
product of $A$ by $B$ with respect to some arbitrary action of $B$ on $A$.

\begin{defn} \label{fjresults} {\rm We define the following notations.

${\cal N}:$ The FICwF$^{{\cal H}^?_*}_{\cal {VC}}$ is true for $\pi_1(M)$ 
for closed nonpositively curved Riemannian manifolds $M$.

${\cal B}:$ The FICwF$^{{\cal H}^?_*}_{\cal {VC}}$ is true for 
${\Bbb Z}^n\rtimes 
{\Bbb Z}$ for all $n$.

${\cal L}:$ If $G=lim_{i\to \infty}G_i$ and the FIC$^{{\cal H}^?_*}_{\cal 
{VC}}$ 
is true for $G_i$ for each $i$ then the FIC$^{{\cal H}^?_*}_{\cal {VC}}$ 
is true for $G$.}\end{defn}

\begin{thm}\label{fj}  $\cal B$, $\cal N$ and $\cal L$ are
satisfied for
the FICwF$^{_P{\cal H}^?_*}_{\cal {VC}}$.\end{thm}
\begin{proof} See theorem 4.8 in \cite{FJ} for
$\cal B$, proposition 2.3 in \cite{FJ} and fact 3.1 in \cite{FR} for
$\cal N$ and theorem 7.1 in \cite{FL} for $\cal L$.\end{proof}


\begin{thm} \label{mth} $(1).$ $\cal B$ implies that the FICwF$^{{\cal 
H}^?_*}_{\cal {VC}}$ is true for the following groups. 

(a). Finitely generated virtually polycyclic groups. 

(b). Finitely generated virtually solvable subgroup of $GL_n({\Bbb C})$. 

In addition if we assume $\cal L$ then we can remove the condition  
`finitely generated' from (a) and (b).

$(2).$ $\cal B$ and $\cal N$ together imply that  
the FICwF$^{{\cal H}^?_*}_{\cal {VC}}$ is true for the following groups.

(a). Cocompact discrete subgroups of virtually connected Lie groups.

(b). Artin full braid groups.

$(3).$ $\cal B$, $\cal N$ and $\cal L$ imply that the FICwF$^{{\cal 
H}^?_*}_{\cal {VC}}$ 
is true for the following groups.

(a). Fundamental groups of Seifert fibered and Haken 
$3$-manifolds. And fundamental groups of arbitrary $3$-manifolds if we 
assume the Geometrization conjecture.

(b). Weak strongly poly-surface groups (see Definition \ref{defn}).

(c). Extensions of closed surface groups by surface groups.

(d). $\pi_1(M)\rtimes \Bbb Z$, where $M$ is a closed Seifert fibered 
space.
\end{thm} 

The following Corollary follows using $(a)$ of Lemma \ref{inverse}.

\begin{cor} Let $G$ be a group containing a finite index subgroup 
$\Gamma$. Assume that $\Gamma$ is isomorphic to a group appearing 
in 
Theorem 
\ref{mth} together with the corresponding hypotheses. Then the 
FICwF$^{{\cal H}^?_*}_{\cal {VC}}$
is true for $G$.\end{cor}

Theorem \ref{fj} and \ref{mth} together imply the following corollary.

\begin{cor} \label{cons} The FICwF$^{_P{\cal H}^?_*}_{\cal {VC}}$ is true 
for 
groups appearing in Theorem \ref{mth}.\end{cor}

Corollary \ref{fj} together with Proposition \ref{corprop} and proposition 
4.10 in \cite{LR} imply the following for the Isomorphism 
Conjecture in the algebraic $K$-theory case.

\begin{cor} \label{kth} Let $G$ be a group which contains a finite 
index subgroup $\Gamma$ so that either $\Gamma$ is a group 
appearing in Theorem \ref{mth} or assume that the FICwF$^{_P{\cal 
H}^?_*}_{\cal {VC}}$ is true for 
$\Gamma$. 
Then for $R={\Bbb Z}$ the algebraic $K$-theory 
assembly map 
$${\cal H}^{G\wr H}_n(E_{\cal {VC}}(G\wr H), {\bf K}_R)\to  
{\cal
H}^{G\wr H}_n(pt, {\bf K}_R)$$ is an isomorphism for $n\leq 1$, where $H$ 
is a finite group. Or equivalently the ICwF$^{_K{\cal   
H}^?_*}_{\cal {VC}}$ is true for
$G$ in dimension $\leq 1$ when $R={\Bbb Z}$.

Here recall   
that the above map is induced by the map $E_{\cal {VC}}(G\wr H)\to pt$ 
and ${\cal 
H}^?_n(X, {\bf K}_R)$ is the standard notation for $_K{\cal 
H}^?_n(X)$.
\end{cor}

\begin{rem}\label{remark}{\rm 
$\cal L$ is known 
for many cases, for example see  theorem 7.1 in \cite{FL} for the 
FIC$^{_K{\cal H}^?_*}_{\cal
{VC}}$ and FIC$^{_L{\cal H}^?_*}_{\cal
{VC}}$ and for the FIC$^{_{KH}{\cal H}^?_*}_{\cal
{VC}}$ see theorem 11.4 in \cite{BL}. We recall now some results 
related to $\cal N$. In \cite{BR}, IC$^{_K{\cal H}^?_*}_{\cal
{VC}}$ is proved for the fundamental groups of closed strictly negatively 
curved Riemannian manifolds. We have already mentioned in the abstract 
that the proof of the FIC$^{_P{\cal H}^?_*}_{\cal {VC}}$ 
is known for all the groups appearing in Theorem \ref{mth} except 
for those in $3(c)$ and $3(d)$  
(the precise references are given during the proof). In case $3(a)$ 
FICwF$^{_P{\cal H}^?_*}_{\cal {VC}}$ is also known. The main ingredients 
behind the proof of Theorem \ref{mth} are Lemmas \ref{begin} to 
\ref{freeproduct}.} \end{rem}

\begin{rem}{\rm We should remark that the FIC$^{_P{\cal H}^?_*}_{\cal 
{VC}}$ was proved for certain class of mapping class groups of surfaces 
of lower genus in 
\cite{BPL}. Assuming $\cal B$, $\cal N$ and $\cal L$ and using the 
methods of this article the 
FICwF$^{{\cal H}^?_*}_{\cal {VC}}$ can also be deduced for these mapping 
class 
groups. We leave the details to the reader.}\end{rem}

\section{Some preliminary results}

\begin{lemma}\label{alglem} [algebraic lemma, \cite{FR}] Let $K$ be a 
finite index normal subgroup of a 
group $G$, then $G$ embeds in $K\wr 
(G/K)$.\end{lemma}

The following lemma is standard.

\begin{lemma} \label{char} (a). Let $G$ be a finite index subgroup of a 
group $K$. Then there is a subgroup $G_1<G$ which is normal and of 
finite index in $K$.

(b). In addition if $K$ is finitely presented then there is a subgroup 
$G_1<G$ which is characteristic and of 
finite index in $K$.\end{lemma}

\begin{lemma} \label{begin} $\cal {B}$ implies that the {\it ${\cal 
P}^{{\cal 
H}^?_*}_{\cal {VC}}$-property} is
satisfied.\end{lemma}

\begin{proof} Note that for two virtually 
cyclic groups $C_1$ and $C_2$, $C_1\times C_2$ contains a free abelian 
(of rank $\leq 2$) normal subgroup (say $A$) of finite index. Hence by 
Lemma 
\ref{alglem} $C_1\times C_2$ is a subgroup of $A\wr F$ for some finite 
group $F$. Therefore $\cal B$ and the hereditary property of the 
FICwF$^{{\cal H}^?_*}_{\cal {VC}}$ completes the 
proof.\end{proof}

\begin{lemma}\label{inverse} Assume that the {\it 
${\cal P}^{{\cal H}^?_*}_{\cal {VC}}$-property} is
satisfied. 

(a). Let $G$ be a finite index subgroup of a group $K$. If the  
FICwF$^{{\cal 
H}^?_*}_{\cal {VC}}$ is true for $G$ then the FICwF$^{{\cal H}^?_*}_{\cal 
{VC}}$ is also true for $K$.

(b). Let $p:G\to Q$ 
be a surjective group homomorphism and assume that the FICwF$^{{\cal 
H}^?_*}_{\cal {VC}}$ is true for $Q$, for ker$(p)$ and for $p^{-1}(C)$ for 
any 
infinite cyclic subgroup $C$ of $Q$. 
Then $G$ 
satisfies the FICwF$^{{\cal H}^?_*}_{\cal {VC}}$.\end{lemma}

\begin{proof} Note that $(a)$ is same as $(2)$ of proposition 5.2 in 
\cite{R4} 
with the only difference that there we assumed $G$ is also normal in $K$. 
But this normality condition can be removed using $(a)$ of Lemma 
\ref{char} and the hereditary property and then applying $(2)$ of 
proposition 5.2 in \cite{R4}.

Also $(b)$ is easily deduced from $(a)$ and $(3)$ of proposition 5.2 in 
\cite{R4}. Since the hypothesis of $(3)$ of proposition 5.2 in 
\cite{R4} was that the FICwF$^{{\cal
H}^?_*}_{\cal {VC}}$ is true for $Q$ and for $p^{-1}(C)$ for 
any virtually 
cyclic subgroup (including the trivial group) $C$ of $Q$. And by 
definition a virtually 
cyclic group is either finite or contains an infinite cyclic subgroup of 
finite index.\end{proof}

An immediate Corollary of $(a)$ is the following.

\begin{cor}\label{corollary} If the ${\cal P}^{{\cal H}^?_*}_{\cal 
{VC}}$-property is satisfied then the FICwF$^{{\cal
H}^?_*}_{\cal {VC}}$ is true for all $H\in {\cal {VC}}$.\end{cor}

\begin{lemma} \label{jsjt} Assume $\cal N$ and that the {\it ${\cal 
P}^{{\cal
H}^?_*}_{\cal {VC}}$-property} is
satisfied. Let $M$ be one of the following 
manifolds. 

(a). A closed Haken 
$3$-manifold such that
there is a hyperbolic piece in the JSJT decomposition of the manifold. 

(b). A compact irreducible $3$-manifold with nonempty
incompressible boundary and there is at least one torus boundary
component. 

(c). A compact surface.

Then the FICwF$^{{\cal H}^?_*}_{\cal {VC}}$ is true for 
$\pi_1(M)$.\end{lemma}

\begin{proof} In $(a)$ by theorem 3.2 in \cite{Le} $M$ supports a 
nonpositively curved metric. The proof is now immediate.

For $(b)$ let $N$ be the double of $M$ along boundary components of 
genus $\geq  2$. Then we can again apply theorem 3.2 in \cite{Le} to see 
that 
the interior of $N$ supports a complete nonpositively curved metric so 
that near the boundary the metric is a product of flat tori and 
$(0,\infty)$. Therefore the double of $N$ is a closed nonpositively curved 
manifold. Now we can apply the hereditary property to complete the proof. 

The proof of $(c)$ also follows by taking the double and using the fact 
that a closed surface either has finite fundamental group or supports a 
nonpositively curved metric.\end{proof}

\begin{lemma} \label{freeproduct} Assume $\cal N$, $\cal L$ and that the 
${\cal P}^{{\cal 
H}^?_*}_{\cal {VC}}$-property is satisfied. Then the followings hold. 
 
(a). 
The  FICwF$^{{\cal
H}^?_*}_{\cal {VC}}$ is true for countable free groups.

(b). 
If the FICwF$^{{\cal 
H}^?_*}_{\cal {VC}}$ is true for two countable groups $G_1$ and $G_2$ then 
the  
FICwF$^{{\cal H}^?_*}_{\cal {VC}}$ is true for the free product 
$G_1*G_2$.\end{lemma}  

We recall here that the above lemma was proved in proposition 5.3 and 
lemma 6.3 in \cite{R4} under a different hypothesis.
 
\begin{proof} For the proof of $(a)$ note that a  
finitely generated free group is isomorphic to the fundamental 
group of a compact surface and hence $(c)$ of Lemma \ref{jsjt} 
applies. And since an infinitely generated countable free group 
is a direct limit of finitely generated free subgroups we can apply $\cal 
L$ in addition to $(c)$ of Lemma \ref{jsjt} to complete the proof.

For the proof of $(b)$ consider the following short exact 
sequence. 
$$1\to K\to G_1*G_2\to 
G_1\times G_2\to 1,$$ where $K$ is the kernel of the homomorphism 
$p:G_1*G_2\to G_1\times G_2$. By $(1)$ of proposition 5.2 in \cite{R4} the 
FICwF$^{{\cal H}^?_*}_{\cal {VC}}$ is true for the product $G_1\times 
G_2$. Now by lemma 5.2 in \cite{R2} $p^{-1}(C)$ is a countable free group 
if $C$ is either the trivial group or an infinite cyclic subgroup of 
$G_1\times G_2$. Therefore we can apply the previous assertion and $(b)$ 
of 
Lemma \ref{inverse} to the above exact sequence to complete the proof. 
\end{proof}

\section{Proof of Theorem \ref{mth}}
At first recall that by Lemma \ref{begin} $\cal B$ implies that the ${\cal 
P}^{{\cal 
H}^?_*}_{\cal {VC}}$-property is satisfied. Therefore we can use Lemma 
\ref{inverse}, Corollary \ref{corollary}, Lemma \ref{jsjt} and Lemma 
\ref{freeproduct} in the proof 
of the Theorem.

\medskip
\noindent
{\it Proof of 1(a).} Let $\Gamma$ be a finitely generated virtually 
polycyclic group. That is, $\Gamma$ contains a (finitely 
generated) polycyclic subgroup of 
finite index. Also a polycyclic group contains a poly-$\Bbb Z$ subgroup of 
finite index (see 5.4.15 in \cite{R}). Therefore by $(a)$ of Lemma 
\ref{inverse} it is enough to prove the 
FICwF$^{{\cal H}^?_*}_{\cal {VC}}$ for poly-$\Bbb Z$ groups. Hence  
we assume that $\Gamma$ is a finitely generated poly-$\Bbb Z$ group. 
Now the 
proof is by induction 
on the virtual cohomological dimension ({\it vcd}) of $\Gamma$.  

We need the following lemma which is easily deduced from 
lemma 4.1 and lemma 4.4 in \cite{FJ}.

\begin{lemma} \label{deduct} Let k=vcd($\Gamma$). Then there is a short 
exact sequence of 
groups. $$1\to \Gamma'\to \Gamma\to A\to 1,$$ where $\Gamma'$ and $A$ 
satisfy the following properties.

(a). $\Gamma'$ is a poly-$\Bbb Z$ group with vcd($\Gamma'$)$\leq k-2$.

(b). There is a short exact sequence. $$1\to A'\to A\to B\to 1,$$ where 
$A'$ contains a finitely generated free abelian subgroup of finite index 
and $B$ is either trivial or infinite cyclic or the infinite dihedral 
group $D_{\infty}$. 
\end{lemma}

Now assume that the FICwF$^{{\cal H}^?_*}_{\cal {VC}}$ is true for 
finitely generated poly-$\Bbb Z$ groups of virtual cohomological 
dimension $\leq k-1$. Then apply the hypothesis and $(b)$ of Lemma 
\ref{inverse} to the two exact sequences in Lemma 
\ref{deduct}. 
The details arguments are routine and we leave it to the reader. This 
completes the proof of $1(a)$.

\medskip
\noindent
{\it Proof of 1(b).} Let $\Gamma$ be a finitely generated virtually 
solvable subgroup of 
$GL_n({\Bbb C})$. Then $\Gamma$ contains a finitely generated torsion 
free
solvable subgroup of finite index. Hence by $(a)$ of Lemma \ref{inverse} 
we can assume that $\Gamma$ is torsion free.  By a theorem of Mal'cev 
( see 15.1.4 in \cite{R}) there 
are subgroups $\Gamma_2<\Gamma_1<\Gamma$ with the following properties.
\begin{itemize}
\item $\Gamma_2$ and $\Gamma_1$ are  
normal in $\Gamma$.
\item $\Gamma_2$ is nilpotent.
\item $\Gamma_1/\Gamma_2$ is abelian.
\item $\Gamma/\Gamma_1$ is finite.
\end{itemize}

Therefore we have the following exact sequences. $$1\to \Gamma_1/\Gamma_2 
\to \Gamma/\Gamma_2 \to \Gamma/\Gamma_1 \to 1.$$ $$1\to \Gamma_2 \to 
\Gamma \to 
\Gamma/\Gamma_2\to 1.$$ Since $\Gamma$ is finitely generated and  
$\Gamma/\Gamma_1$ is finite, it follows that $\Gamma_1$ is also finitely 
generated. Hence $\Gamma_1/\Gamma_2$ is finitely generated abelian. 
Consequently $\Gamma/\Gamma_2$ contains a finitely generated free abelian 
subgroup of finite index. Since $\Gamma_2$ is torsion free nilpotent it is 
poly-$\Bbb Z$. Now using $1(a)$, $\cal B$ and $(b)$ of Lemma \ref{inverse} 
we complete the proof. 

\medskip
\noindent
{\it Proof of 2(a).} We follow the structure of the proof of theorem 2.1 
in \cite{FJ}. So 
let $\Gamma$ be a discrete cocompact subgroup of a virtually connected Lie 
group $G$. First let us show that we can assume that the Lie group has no 
nontrivial compact connected normal subgroup. So let $C$ denotes the 
maximal compact connected normal subgroup of $G$. Let $\Gamma'=q(\Gamma)$ 
and $F=\Gamma\cap C$ where $q:G\to G/C$ is the quotient map. Then $F$ is a 
finite normal subgroup of $\Gamma$ with quotient $\Gamma'$ which is a 
discrete cocompact subgroup of the virtually connected Lie group $G/C$. 
Note that $G/C$ contains no nontrivial compact connected normal subgroup. 
Applying $(b)$ of Lemma \ref{inverse} to the map $q$ we see that the 
FICwF$^{{\cal H}^?_*}_{\cal {VC}}$ for $\Gamma'$ implies the FICwF$^{{\cal 
H}^?_*}_{\cal {VC}}$ for $\Gamma$. Since if a group contains a finite 
normal subgroup with infinite cyclic quotient then the group is virtually 
infinite cyclic.

Therefore from now on we assume that $G$ has no nontrivial compact 
connected normal subgroup. Next we deduce that we can also assume that the 
identity component $G_e$ of $G$ is a semisimple Lie group. By the Levi 
decomposition we have that $G$ contains a maximal closed connected normal 
solvable subgroup $R$ with quotient $S$ which is a virtually connected Lie 
group with the identity component $S_e$ semisimple. That is we have the 
following exact sequence. $$1\to R\to G\to S\to 1.$$ Let 
$\Gamma_R=\Gamma\cap R$ and $\Gamma_S=p(\Gamma)$, where $p$ denotes the 
homomorphism $G\to S$. Then we have the 
following (see 2.6 (a), (b) and (c) in \cite{FJ}). A short exact 
sequence $$1\to \Gamma_R\to \Gamma\to \Gamma_S\to 1,$$ where $\Gamma_R$ is 
virtually poly-$\Bbb Z$ and discrete cocompact in $R$ and $\Gamma_S$ is  
discrete cocompact in  
$S$. Therefore using case $1(a)$ of the Theorem and by $(b)$ of Lemma 
\ref{inverse} 
we see that it is enough to prove the FICwF$^{{\cal H}^?_*}_{\cal {VC}}$ 
for $\Gamma_S$. 

Henceforth we assume that $G$ has no nontrivial compact
connected normal subgroup and $G_e$ is semisimple. Let $H=\{g\in G\ |\ 
gh=hg$ for all $h\in G_e\}$ and $q:G\to G/H$ be the quotient map. Then we 
have the following (see 2.5 (a), (b), (c) and (d) in \cite{FJ}). $$1\to 
\Gamma_H=H\cap\Gamma\to 
\Gamma\to \Gamma'=q(\Gamma)\to 1,$$ where $\Gamma_H$ is a finite extension 
of a finitely generated abelian group, $\Gamma'$ is a cocompact discrete 
subgroup of $G/H$ and $G/H$ is a virtually connected linear Lie group such 
that the identity component of $G/H$ is semisimple. Now note that 
$\Gamma'$ acts cocompactly and properly discontinuously on a 
complete nonpositively curved Riemannian manifold (namely on 
$(G/H)/C$, where $C$ is maximal compact subgroup of $G/H$). Hence by 
$\cal N$ the FICwF$^{{\cal H}^?_*}_{\cal {VC}}$ is true for $\Gamma'$. 
Since 
$\Gamma_H$ contains a finitely generated free abelian subgroup of finite 
index, we can therefore apply $1(a)$ of the Theorem and Lemma 
\ref{inverse} to the above exact sequence to complete the proof 
of $2(a)$. 

\medskip
\noindent
{\it Proof of 2(b).} Let 
$$M_n=\{(x_1,\ldots ,x_{n+1})\in {\Bbb C}^{n+1}\ |\ x_i\neq x_j,\ 
\text{for}\ i\neq 
j\}.$$ Then the symmetric group $S_{n+1}$ acts on $M_n$ freely. By 
definition the Artin full braid group (denoted by $B_n$) on $n$ strings is 
the fundamental group of the manifold $M_n/S_{n+1}$ and the pure braid 
group is defined as the fundamental group of $M_n$ and is denoted by 
$PB_n$. Hence $PB_n$ is a normal subgroup of $B_n$ with quotient 
$S_{n+1}$. Using $(a)$ of Lemma \ref{inverse} we see that it is enough to 
prove the 
FICwF$^{{\cal H}^?_*}_{\cal {VC}}$ for $PB_n$. The proof is by induction 
on $n$. For $n=1$, $PB_n\simeq {\Bbb Z}$ and hence by Corollary
\ref{corollary} we can start the 
induction. So assume the result for $PB_k$ for $k\leq n-1$. We will show 
that the FICwF$^{{\cal H}^?_*}_{\cal {VC}}$ is true for $PB_n$. Consider 
the following projection map. $$p:M_n\to M_{n-1}$$$$(x_1,\ldots , 
x_{n+1})\mapsto (x_1,\ldots ,
x_n).$$ This map is a locally trivial fiber bundle projection with fiber 
diffeomorphic to $N={\Bbb C}-\{$$n$-points$\}$. Therefore we have the 
following exact sequence. $$1\to 
\pi_1(N)\to PB_n\to PB_{n-1}\to 1.$$ The following is easy to show.
\begin{itemize}
\item For an infinite cyclic subgroup $C$ of $PB_{n-1}$, 
$p_*^{-1}(C)$ is isomorphic  
to $\pi_1(N_f)$ where $N_f$ denotes the 
mapping torus of the monodromy diffeomorphism $f:N\to N$ corresponding to 
a generator of $C$. It follows that $N_f$ is diffeomorphic to the interior 
of a compact Haken $3$-manifold with incompressible boundary 
components. 
\end{itemize}

Note that all the boundary components of $N_f$ could be Klein bottles. In 
that 
case we can take a finite sheeted cover of $N_f$ which is a Haken 
$3$-manifold with incompressible tori boundary components.

Now we can apply $(a)$ and $(b)$ of Lemma \ref{inverse}, $(b)$ of Lemma 
\ref{jsjt} 
and $(a)$ of Lemma \ref{freeproduct}  
to 
complete the proof of $2(b)$.

\medskip
\noindent
{\it Proof of 3(a).} The 
proof goes in the same line as the proof of the FICwF$^{_P{\cal 
H}^?_*}_{\cal {VC}}$ for $3$-manifold groups (see \cite{R2} and 
\cite{R3}). 

First of all using $\cal L$ we can assume that the $3$-manifold is 
compact. Also applying $(a)$ of Lemma \ref{inverse} we can 
assume that the $3$-manifold is orientable and that all boundary 
components are orientable. So let $M$ be a compact 
orientable $3$-manifold with orientable boundary. If there is any sphere 
boundary then we can cap these boundaries by $3$-discs. This does not 
change the 
fundamental group. Therefore we can assume that either $M$ is closed or 
all boundary components are orientable surfaces of genus $\geq 1$. By 
lemma 3.1 in \cite{R2}  $\pi_1(M)\simeq \pi_1(M_1)*\cdots 
*\pi_1(M_k)*F^r$, where $M_i$ is a compact orientable irreducible 
$3$-manifold for each $i=1,2,\ldots , k$ and $F^r$ is a free group of rank 
$r$. Hence applying Lemma \ref{freeproduct} we are reduced to the 
situation of compact orientable 
irreducible $3$-manifold $M$. Also we can assume $\pi_1(M)$ is infinite 
as the FICwF$^{{\cal H}^?_*}_{\cal {VC}}$ is true for finite groups. 

Now by Thurston's Geometrization conjecture $M$ has the following 
possibilities.

\begin{itemize}
\item Closed Seifert fibered space.
\item Haken manifold.
\item Closed hyperbolic manifold.
\end{itemize} 

Using $\cal N$ we only have to consider the first two cases. In the first 
case there is an exact sequence. $$1\to {\Bbb Z}\to \pi_1(M)\to 
\pi_1^{orb}(S)\to 1,$$ where $S$ is the base orbifold of the Seifert 
fibered space $M$. It is well known that $\pi_1^{orb}(S)$ is 
either finite or contains a closed  
surface subgroup of finite index. Now we can apply $(c)$ of Lemma 
\ref{jsjt} and $(a)$ and $(b)$ of 
Lemma \ref{inverse} to complete the proof in this case. 

For the Haken manifold situation let us first consider the nonempty 
boundary case. By lemma 6.4 in \cite{R2} we have that $\pi_1(M)\simeq 
\pi_1(M_1)*\cdots *\pi_1(M_l)*F^s$, where each $M_i$ is compact orientable 
irreducible and has incompressible boundary and $F^s$ is a free group of 
rank $s$. We apply $(a)$ and $(b)$ of Lemma \ref{freeproduct}, $(b)$ of 
Lemma \ref{jsjt},  
theorem 1.1.1 in \cite{R3}, $\cal N$ and the hereditary property of the 
FICwF$^{{\cal H}^?_*}_{\cal {VC}}$  to complete the proof in this case. 
(Recall here that the theorem 1.1.1 and the remark 1.1.1 in 
\cite{R3} say that the fundamental 
group of each $M_i$ is a subgroup of the fundamental group of a closed 
$3$-manifold $P$ which is either Seifert fibered or supports 
a nonpositively curved metric.)

Finally we assume that $M$ is a closed Haken manifold. If $M$ contains no 
incompressible torus then by Thurston's hyperbolization theorem $M$ is 
hyperbolic. Hence using $\cal N$ and $(a)$ of Lemma 
\ref{jsjt} we can assume that $M$ contains an 
incompressible torus and has no hyperbolic piece in the JSJT 
decomposition. Therefore $M$ is a graph manifold. By theorem 1.1 in 
\cite{L} either $M$ has a finite sheeted cover which is a torus bundle 
over the circle or for any positive integer $k$ there is a finite sheeted 
cover $M_k$ of $M$ so that rank of $H_1(M_k, {\Bbb Z})$ is $\geq k$. In 
the first possibility we can apply $1(a)$ of the Theorem to 
complete the proof. For the second case let $k\geq 2$ and consider the
following exact 
sequence. $$1\to [\pi_1(M_k),\pi_1(M_k)]\to \pi_1(M_k)\to H_1(M_k, {\Bbb 
Z})\to 1.$$ By $(a)$ of Lemma \ref{inverse} it is enough to 
prove the FICwF$^{{\cal H}^?_*}_{\cal {VC}}$ for $\pi_1(M_k)$. We will 
apply $(b)$ of Lemma \ref{inverse} to the above exact sequence. So Let $C$ 
be either the trivial group or an infinite  
cyclic subgroup of $H_1(M_k, {\Bbb Z})$. 
Then $A^{-1}(C)\simeq \pi_1(K)$ where $K$ is a noncompact irreducible 
(being an infinite sheeted covering of an irreducible $3$-manifold) 
$3$-manifold. Here $A$ 
denotes the homomorphism $\pi_1(M_k)\to H_1(M_k, {\Bbb
Z})$. Let 
$K=\cup 
_iK_i$ where $K_i$'s are increasing union of non-simply connected compact 
submanifold of $K$. By lemma 6.3 and lemma 6.4 in \cite{R2} we have that 
$\pi_1(K_i)\simeq
\pi_1(M_1)*\cdots *\pi_1(M_l)*F^s$, where each $M_i$ is compact orientable
irreducible and has incompressible boundary and $F^s$ is a free group of
rank $s$. Therefore we can use $\cal L$, Lemma \ref{freeproduct} and Haken 
manifold with nonempty 
boundary case as in the previous paragraph to complete the proof of 
$2(b)$. 

\medskip
\noindent
{\it Proof of 3(b).}
Let us recall the definition of weak strongly poly-surface group. 

\begin{defn} \label{defn} [definition 1.2.1 in \cite{R3}] {\rm A discrete 
group $\Gamma$ is called 
{\it weak strongly poly-surface} if there exists a finite filtration of
$\Gamma$ by subgroups: $1=\Gamma_0\subset \Gamma_1\subset \cdots \subset
\Gamma_n=\Gamma$ such that the following conditions are satisfied:

\begin{itemize}
\item $\Gamma_i$ is normal in $\Gamma$ for each $i$.

\item $\Gamma_{i+1}/\Gamma_i$ is isomorphic to the fundamental group
of a surface $F_i$ (say).

\item for each $\gamma\in \Gamma$ and $i$  there is a diffeomorphism
$f:F_i\to F_i$ such
that the induced automorphism $f_{\#}$ of $\pi_1(F_i)$ is equal to
$c_\gamma$ up to inner automorphism, where $c_\gamma$ is the 
automorphism
of $\Gamma_{i+1}/\Gamma_i\ \simeq \ \pi_1(F_i)$ induced by the conjugation
action on $\Gamma$ by $\gamma$.
\end{itemize}

In such a situation we say that the group $\Gamma$ has {\it rank} $\leq
n$}.\end{defn}

The proof is by induction on the rank $n$ of the weak strongly 
poly-surface group $\Gamma$. Therefore assume that the FICwF$^{{\cal 
H}^?_*}_{\cal {VC}}$ is true for all weak strongly
poly-surface group of rank $\leq n-1$ and let $\Gamma$ has rank $n$. 
Consider the following exact sequence. $$1\to \Gamma_1\to \Gamma\to 
\Gamma/\Gamma_1\to 1.$$ The followings are easy to verify. 
\begin{itemize}
\item $\Gamma/\Gamma_1$ is weak strongly poly-surface and has rank $\leq 
n-1$. 
\item $q^{-1}(C)$ is either a surface group or isomorphic to the 
fundamental group of a $3$ manifold, where $C$ is either the trivial 
group or an infinite cyclic subgroup of 
$\Gamma/\Gamma_1$ respectively.
\end{itemize}

Now we can apply the induction hypothesis, $(b)$ of Lemma \ref{inverse}, 
$(c)$ of Lemma \ref{jsjt} and $3(a)$ of the Theorem to 
show that the  FICwF$^{{\cal
H}^?_*}_{\cal {VC}}$ is true for $\Gamma$.

\medskip
\noindent
{\it Proof of 3(c).}
Let $\Gamma$ be an extension of a closed surface group by a surface 
group. Hence we have the following. $$1\to \pi_1(S)\to \Gamma\to 
\pi_1(S')\to 1$$ where $S$ is a closed surface and $S'$ is a surface. 

Using $(a)$ of Lemma \ref{inverse} and $(c)$ of Lemma \ref{jsjt} we can 
assume that $\pi_1(S')$ is not finite 
(and hence torsion free). Also since extension of finite group by 
infinite cyclic group is virtually cyclic we can also assume that 
$\pi_1(S)$ is infinite (and hence again torsion free). Recall that an 
automorphism of a closed surface group is induced by a diffeomorphism of 
the surface. Hence for any infinite cyclic subgroup $C$ of $\pi_1(S')$,  
$p^{-1}(C)$ is isomorphic to the fundamental group of a $3$-manifold. 
Therefore we can apply $(b)$ of Lemma \ref{inverse} and $3(a)$ of the 
Theorem to complete the proof of $3(c)$.

\medskip
\noindent
{\it Proof of 3(d).}
Let $\Gamma=\pi_1(M)\rtimes \Bbb Z$, where $M$ is 
a closed Seifert fibered space. Using $(b)$ of Lemma \ref{char} and $(a)$ 
of Lemma \ref{inverse} we can 
assume that $M$ and its base surface $S$ are both orientable. This implies 
that the cyclic subgroup (say $C$) of $\pi_1(M)$ generated by a regular 
fiber is central (see \cite{H}). Now if $\pi_1^{orb}(S)$ is finite then 
$\Gamma$ is virtually poly-$\Bbb Z$. Therefore using case $1(a)$ of the 
Theorem we can assume that $\pi_1^{orb}(S)$ is infinite and hence $C$ is 
also infinite cyclic. Since $C$ is central we have the following exact 
sequence. $$1\to 
{\Bbb Z}\to \pi_1(M)\rtimes {\Bbb Z}\to \pi_1^{orb}(S)\rtimes {\Bbb Z}\to 
1.$$ 

As $\pi_1^{orb}(S)$ is infinite, it contains a closed 
surface subgroup of finite index. Now using $(b)$ of Lemma \ref{char} it 
is easy to find a characteristic closed
surface subgroup of  $\pi_1^{orb}(S)$ of finite index. Therefore 
$\pi_1^{orb}(S)\rtimes 
{\Bbb 
Z}$ contains a $3$-manifold group of finite index. Now we can apply 
$(a)$ and $(b)$ of Lemma 
\ref{inverse} and $1(a)$ and $3(a)$ of the Theorem to complete the proof 
of $3(d)$.

This completes the proof of Theorem \ref{mth}. 

\newpage
\bibliographystyle{plain}
\ifx\undefined\bysame
\newcommand{\bysame}{\leavevmode\hbox to3em{\hrulefill}\,}
\fi

\medskip

\end{document}